\documentclass[a4paper,11pt,leqno,twoside]{article}
\usepackage{amsmath,amsthm,amssymb}
\usepackage[usenames]{color}
\usepackage[all]{xy}
\numberwithin{equation}{section}

\newcommand{\bC}{\mathbb{C}}

\newcommand{\bN}{\mathbb{N}}

\newcommand{\bQ}{\mathbb{Q}}

\newcommand{\bZ}{\mathbb{Z}}

\newcommand{\fS}{\mathfrak{S}}

\renewcommand{\a}{\alpha}
\renewcommand{\b}{\beta}

\renewcommand{\d}{\delta}
\newcommand{\e}{\varepsilon}
\newcommand{\z}{\zeta}

\renewcommand{\k}{\kappa}
\renewcommand{\l}{\lambda}
\newcommand{\m}{\mu}

\newcommand{\s}{\sigma}
\renewcommand{\t}{\tau}
\renewcommand{\c}{\chi}
\newcommand{\vp}{\varphi}
\newcommand{\om}{\omega}

\newcommand{\G}{\Gamma}

\renewcommand{\Re}{\mathrm{Re}\,}

\newcommand{\Gauss}[1]{\lfloor{#1}\rfloor}

\newcommand{\boldtitle}[1]{\title{\bfseries #1}}

\newenvironment{MSC}{%
\smallbreak
\noindent \textbf{2000\ Mathematics Subject Classification\,:}\ 
}

\newenvironment{keywords}{%
\noindent\textbf{Key words and phrases\,:}\itshape}

\newenvironment{Acknowledgement}{%
\noindent\textit{Acknowledgement.}}

\theoremstyle{theorem}
\newtheorem*{multitheorem}{\variable@name}

\theoremstyle{definition}
\newtheorem*{multiproclaim}{\variable@name}

\theoremstyle{plain}
\newtheorem{thm}{Theorem}[section]
\newtheorem{prop}[thm]{Proposition}
\newtheorem{lem}[thm]{Lemma}
\newtheorem{cor}[thm]{Corollary}

\newtheorem{ques}[thm]{Question}

\theoremstyle{definition}

\newtheorem{example}[thm]{Example}
\newtheorem{remark}[thm]{Remark}

\newcommand{\bsym}{\boldsymbol}

\newcommand{\ul}[1]{\underline{#1}}
\newcommand{\sgn}{\mathrm{sgn}}

\newcommand{\ol}[1]{\overline{#1}}
\renewcommand{\pmod}[1]{(\mathrm{mod}\ {#1})}
\newcommand{\wt}[1]{\widetilde{#1}}

\newcommand{\bul}{\bullet}
\newcommand{\odd}{\mathrm{o}}
\newcommand{\even}{\mathrm{e}}
\newcommand{\tri}{\mathrm{tri\,}}
\newcommand{\cosec}{\mathrm{cosec\,}}

\newcommand{\DS}[1]{\displaystyle{#1}}

\boldtitle{Evaluations of multiple Dirichlet $L$-values\\
via symmetric functions}
\author{Yoshinori Yamasaki
\footnote
{Partially supported by Grant-in-Aid for JSPS Fellows No.\,19002485.}} 
\date{\today}

\begin{document}

\setlength{\baselineskip}{15pt}
\maketitle

\begin{abstract}
 We explicitly evaluate a special type of multiple Dirichlet $L$-values
 at positive integers in two different ways:
 One approach involves using symmetric functions,
 while the other involves using a generating function of the values.
 Equating these two expressions,
 we derive several summation formulae
 involving the Bernoulli and Euler numbers.
 Moreover,
 values at non-positive integers, called central limit values,
 are also studied.   
\begin{MSC}
 {\it Primary} 11M41;   
 {\it Secondary} 11B68. 
\end{MSC} 
\begin{keywords}
 multiple $L$-values,
 multiple $q$-$L$-values,
 central limit values,
 Dirichlet $L$-functions,
 Bernoulli numbers,
 Euler numbers,
 symmetric functions.
\end{keywords}
\end{abstract}

\section{Introduction}
\label{sec:introduction}

%
%

 Let $d\in\bN$,
 $s_1,\ldots,s_d\in\bC$
 and $f_1,\ldots,f_d$ be certain arithmetic functions such as Dirichlet characters.
 In this paper,
 for each $\om\in\{\bul,\star\}$,
 we study the multiple $L$-function      
\[
 L^{\om}_d(s_1,\ldots,s_d;f_1,\ldots,f_d)
:=\sum_{(m_{1},\ldots,m_{d})\in I^{\om}_{d}}
\frac{f_{1}(m_{1})\cdots f_{d}(m_{d})}{m_{1}^{s_{1}}\cdots m_{d}^{s_{d}}},
\]
 where $I^{\om}_{d}:=\lim_{M\to +\infty}I^{\om}_d(M)$ with 
\[
 I^{\om}_d(M):=
\begin{cases}
 \bigl\{(m_1,\ldots,m_d)\in \bN^{d}\,|\,m_1<\cdots<m_d\le M\bigr\} &
 \ \ \textrm{($\om=\bul$)},\\
 \bigl\{(m_1,\ldots,m_d)\in \bN^{d}\,|\,m_1\le\cdots\le m_d\le M\bigr\} &
 \ \ \textrm{($\om=\star$)}.
\end{cases}
\]
 In particular,
 when $f_j=\bsym{1}$ for all $j$
 where $\bsym{1}$ is the trivial character,
 we write $L^{\om}_d$ as $\z^{\om}_d$
 and call it the multiple zeta function
 ($\z^{\star}_d$ is often called
 the multiple zeta-star function
 or the non-strict multiple zeta function).
 Moreover, when $d=1$, we write $L(s;f)=L^{\om}_1(s;f)$,
 whence $L(s;\bsym{1})=\z(s)$ where $\z(s)$ is the Riemann zeta function.  
 
 Let $\{x\}^{d}$ be the $d$-times copy of the variable $x$.
 Though it is, in general, difficult to evaluate the value
 $\z^{\om}_d(k_1,\ldots,k_d)$ for given positive integers $k_1,\ldots,k_d$,
 there are many explicit results
 concerning the values $\z^{\om}_d(\{2k\}^d)$
 (see 
 \cite{BorweinBradleyBroadhurst1997,BowmanBradley2001,Kamano2006,Nakamura2005,Wenchang2000,Zagier1994}  
 for the case $\om=\bul$ and \cite{Muneta,Ohno,Ohno2} for $\om=\star$).
 Let $\c$ be a Dirichlet character and
 $\k$ an integer with the same parity as that of $\c$.
 The purpose of the present paper is,
 as a generalization and unification of the
 above results,
 to obtain explicit descriptions of the values
 $L^{\om}_d(\{\k\}^d;\{\c\}^d)$
 from a symmetric functions' viewpoint.
 That is, we evaluate the values
 $L^{\bul}_d(\{\k\}^d;\{\c\}^d)$ and $L^{\star}_d(\{\k\}^d;\{\c\}^d)$
 as specializations of the elementary and complete symmetric function,
 respectively.  

 When $\k$ is positive,
 we do this in the following two distinct ways:
 In the first method, 
 making use of the relations among the above symmetric functions and
 the power-sum symmetric function,
 we express the values $L^{\om}_d(\{\k\}^d;\{\c\}^d)$
 in terms of a sum over partitions of $d$ (Theorem~\ref{thm:evaluation1}).
 In the second method,
 we consider a generating function of $L^{\om}_d(\{\k\}^d;\{\c\}^d)$. 
 From the product expressions of the generating function,
 we show that $L^{\om}_d(\{\k\}^d;\{\c\}^d)$ can be expressed
 in terms of $A^{\om}_{\k d}(\k;\c)$, 
 which is defined by a finite sum
 involving the multinomial coefficients and roots of unity
 (Theorem~\ref{thm:evaluation2}).
 As an application,
 in Section~\ref{sec:applications},
 equating the two expressions,
 we give several summation formulae involving
 the Bernoulli numbers $B_{n}$ and the Euler numbers $E_{n}$.
 It should be noted that some of these
 (formulae \eqref{for:sum-Bernoulli1}, \eqref{for:sum-Bernoulli2} and
 \eqref{for:sum-Euler})
 can be regarded as solutions of a kind of ``inverse problem'',
 as previously studied in \cite{Haukkanen1994}
 (see also \cite{Dilcher1996}),
 for the sequences $B_{2k}$ and $E_{2k}$. 

 Moreover,
 we also study a special value, called a {\it central limit value}
 (see \cite{AkiyamaTanigawa2001,AkiyamaEgamiTanigawa2001}),
 of the function $L^{\om}_d(s_1,\ldots,s_d;\{\c\}^d)$
 at $s_1=\cdots=s_d=-\k$ where $\k$ is a non-negative integer.
 Note that such a value can be defined 
 because $L^{\om}_d(s_1,\ldots,s_d;\{\c\}^d)$ admits a meromorphic
 continuation to the whole space $\bC^d$.
 We then obtain formulae  that extend
 the well-known expressions concerning the value
 $L(-\k,\c)$ of the usual Dirichlet $L$-function
 (Theorem~\ref{thm:central values}).    

 In the final section,
 we discuss some generalizations of $L^{\om}_d$,
 and also give questions 
 for the future study.

 Note that,
 for completeness and the reader's convenience,
 we restate or briefly prove the corresponding results for $\c=\bsym{1}$. 

 Throughout the present paper,
 we use the following notation for Dirichlet characters:
 Let $\c$ be a Dirichlet character modulo $N$.
 Denote by
 $\c'$ the primitive character that induces $\c$,
 $N(\c)$ the conductor,
 $e(\c)\in\{0,1\}$ the parity, i.e., $\c(-1)=(-1)^{e(\c)}$
 (we call $\c$ even  if $e(\c)=0$ and odd otherwise),
 $\t(\c):=\sum^{N}_{a=1}\c(a)e^{\frac{2\pi ia}{N}}$ the Gauss sum of $\c$ and
 $B_{n,\c}$ the generalized Bernoulli number associated with $\c$
 defined by 
 $\sum^{N}_{a=1}\frac{\c(a)te^{at}}{e^{Nt}-1}=\sum^{\infty}_{n=0}B_{n,\c}\frac{t^n}{n!}$.  
 Note that
 $B_{n,\bsym{1}}=B_n$ and $B_{n,\c_{-4}}=-\frac{1}{2}nE_{n-1}$
 where $\c_{-4}$ is the primitive Dirichlet character modulo $4$.
 Further, for simplification,
 we always write $\k=\k(k,\c):=2k+e(\c)$ with $k\in\bZ$.

\section{Multiple $L$-values via the symmetric functions}
\label{sec:symmetric functions}

\subsection{Symmetric functions}

 Let $\l=(\l_1,\l_2,\ldots)$ be a partition.
 We write $\l\vdash d$ if $\l$ is a partition of $d$,
 that is, $|\l|:=\sum_{j\ge 1}\l_j=d$ and
 put $\ell(\l):=\max\{j\,|\,\l_j\ne 0\}$ (the length of $\l$) and
 $m_i(\l):=\#\{j\,|\l_j=i\,\}$ (the multiplicity of $i$ in $\l$)
 as usual. 
 We denote by $\l_{\odd}$ (resp. $\l_{\even}$) the partition
 whose parts consist of all the odd (resp. even) ones of $\l$
 and set
 $\e_{\l}:=(-1)^{|\l|-\ell(\l)}$ and
 $z_{\l}:=\prod_{j\ge 1}j^{m_j(\l)}m_j(\l)!$.

 For a partition $\l=(\l_1,\l_2,\ldots)$,
 the elementary symmetric function $e_{\l}$,
 the complete symmetric function $h_{\l}$ and
 the power-sum symmetric function $p_{\l}$
 in $\bsym{x}=(x_1,x_2,\ldots)$
 are respectively defined by
 $e_{\l}(\bsym{x}):=\prod_{j\ge 1}e_{\l_j}(\bsym{x})$,
 $h_{\l}(\bsym{x}):=\prod_{j\ge 1}h_{\l_j}(\bsym{x})$ and
 $p_{\l}(\bsym{x}):=\prod_{j\ge 1}p_{\l_j}(\bsym{x})$ where,
 for $r\in\bN$,
\[
 e_{r}:=\sum_{n_1<\cdots<n_r}x_{n_1}\cdots x_{n_r}, \quad 
 h_{r}:=\sum_{n_1\le \cdots\le n_r}x_{n_1}\cdots x_{n_r}
 \quad \textrm{and} \quad 
 p_{r}:=\sum^{\infty}_{n=1} x_n^r.
\]
 Put $v^{\bul}_{\l}:=e_{\l}$, $v^{\star}_{\l}:=h_{\l}$,
 $\e^{\bul}_{\l}:=\e_{\l}$, $\e^{\star}_{\l}:=1$,
 $\e_{\bul}:=-1$ and $\e_{\star}:=1$.
 Then, for each $\om\in\{\bul,\star\}$, we have
\begin{equation}
\label{for:ephp}
 v^{\om}_{(d)}
=\sum_{\m\vdash d}\frac{\e^{\om}_{\m}}{z_{\m}}p_{\m}.
\end{equation}
 This follows from the following identities of
 the generating function of $v^{\om}_{(d)}(\bsym{x})$
 (for further details,
 see the proof of $(2.14)$ and $(2.14)'$
 in \cite[p.\,$25$]{Macdonald1995}); 
\begin{equation}
\label{for:gen}
\sum^{\infty}_{d=0}\e_{\om}^{d}v^{\om}_{(d)}t^d
=\prod^{\infty}_{n=1}(1-x_nt)^{-\e_{\om}}
=\exp\Bigl(\e_{\om}\sum^{\infty}_{n=1}\frac{1}{n}p_nt^n\Bigr)
=\sum^{\infty}_{d=0}\e_{\om}^d
\Bigl\{\sum_{\m\vdash d}\frac{\e^{\om}_{\m}}{z_{\m}}p_{\m}\Bigr\}t^d.
\end{equation}

\subsection{Evaluation formula $\mathrm{I}$}

 For a partition $\m=(\m_1,\m_2,\ldots)$
 and $a,b\in\bZ$ satisfying $a\m_{\ell(\m)}+b\ge0$,
 let $\c^{j}_{\m}:=(\c^{\m_j})'$ for $1\le j\le \ell(\m)$ and set 
\[
 N_{\m}(\c)
:=\prod^{\ell(\m)}_{j=1}N(\c^{j}_{\m})^{\m_j}, \ \ 
 e_{\m}(\c)
:=\sum^{\ell(\m)}_{j=1}e(\c^{j}_{\m}), \ \  
 \t_{\m}(\c)
:=\prod^{\ell(\m)}_{j=1}\t(\c^{j}_{\m}) \ \ 
\textrm{and} \ \   
 \wt{B}_{a\m+b,\c}
:=\prod^{\ell(\m)}_{j=1}\frac{B_{a\m_j+b,\c^{j}_{\m}}}{(a\m_j+b)!}.
\]
 We similarly put
 $\wt{B}_{a\m+b}:=\prod^{\ell(\m)}_{j=1}\frac{B_{a\m_j+b}}{(a\m_j+b)!}$ 
 and 
 $\wt{E}_{a\m+b}:=\prod^{\ell(\m)}_{j=1}\frac{E_{a\m_j+b}}{(a\m_j+b)!}$.
 Moreover, set 
\[
 \a_{\c}(s)
:=\prod_{p:\textrm{prime}\atop p|N,p\nmid N(\c')}\bigl(1-\c'(p)p^{-s}\bigr)
 \ \ \textrm{and} \ \    
 W_{\m}(s;\c)
:=\prod^{\ell(\m)}_{j=1}\a_{\c^{\m_j}}(s\m_j).
\]
 Using these definitions, we obtain the following theorem: 

\begin{thm}
\label{thm:evaluation1}
 Let $\c$ be a Dirichlet character modulo $N$ and
 $\k=2k+e(\c)\ge 1$. Then, we have  
\begin{equation}
\label{for:multiple Dirichlet}
 L^{\om}_{d}\bigl(\{\k\}^d;\{\c\}^d\bigr)
=(-1)^{\frac{\k d}{2}}2^{\k d}
\Biggl\{
\sum_{\m\vdash d}\frac{\e^{\om}_{\m}}{z_{\m}}
\frac{(-1)^{\ell(\m)-e_{\m}(\c)}}{2^{\ell(\m)}}
\frac{W_{\m}(\k;\c)\t_{\m}(\c)}{N_{\m}(\c)^{\k}}\wt{B}_{\k\m,\ol{\c}}
\Biggr\}\pi^{\k d}.
\end{equation}
\end{thm}
\begin{proof}
 We first recall the case $d=1$ (see, e.g., \cite{Neukirch1999}).
 If $\c$ is primitive, then we have 
\begin{equation}
\label{for:primitive}
 L\bigl(\k;\c\bigr)
=(-1)^{k+1-\frac{e(\c)}{2}}\frac{\t(\c)}{2}\Bigl(\frac{2\pi}{N(\c)}\Bigr)^{\k}
\frac{B_{\k,\ol{\c}}}{\k!}.
\end{equation}
 Otherwise, from the Euler product expression,
 we have $L\bigl(\k;\c\bigr)=\a_{\c}(\k)L\bigl(\k;\c'\bigr)$.
 Hence, in this case,
 the claim \eqref{for:multiple Dirichlet} clearly holds. 

 Now, let $d\ge 2$.
 Specializing $x_n=\c(n)n^{-s}$ with $\Re(s)>1$, we have  
\begin{align}
\label{for:L-symmetric}
 v^{\om}_{(d)}(\bsym{x})
=L^{\om}_{d}\bigl(\{s\}^d;\{\c\}^d\bigr)
\quad \textrm{and} \quad 
 p_{(d)}(\bsym{x})
=L(ds;\c^d).
\end{align}
 Note that
 these series also converge for $s=1$
 if $\c$ is not principal (see \cite{ArakawaKaneko2004}). 
 Then, from expression \eqref{for:ephp}, we have
\begin{equation} 
\label{for:L}
 L^{\om}_{d}\bigl(\{s\}^d;\{\c\}^d\bigr)
=\sum_{\m\vdash d}
\frac{\e^{\om}_{\m}}{z_{\m}}\prod^{\ell(\m)}_{j=1}L(\m_js;\c^{\m_j})
=\sum_{\m\vdash d}
\frac{\e^{\om}_{\m}}{z_{\m}}W_{\m}(s,\c)
\prod^{\ell(\m)}_{j=1}L(\m_js;\c^{j}_{\m}).
\end{equation}
 We here put $s=\k=\k(k,\c)\ge 1$.
 Notice that $e(\c^{j}_{\m})\equiv e(\c)\m_j$ $\pmod{2}$ because 
 $(-1)^{e(\c^{j}_{\m})}=\c^{j}_{\m}(-1)=\c(-1)^{\m_j}=(-1)^{e(\c)\m_j}$, 
 and hence  
\begin{equation}
\label{for:kappa}
 \k\m_j=\k(k,\c)\m_j=
\begin{cases}
 \k(k\m_j,\c^{j}_{\m}) & \textrm{if $e(\c)=0$},\\
 \k\bigl(k\m_j+\frac{\m_j}{2},\c^{j}_{\m}\bigr) &
 \textrm{if $e(\c)=1$ and $\m_j$ is even},\\ 
 \k\bigl(k\m_j+\frac{\m_j-1}{2},\c^{j}_{\m}\bigr) &
 \textrm{if $e(\c)=1$ and $\m_j$ is odd}.
\end{cases}
\end{equation}
 Therefore,
 from equations \eqref{for:primitive} and \eqref{for:kappa},
 we have
\begin{equation}
\label{for:L-middd}
 L\bigl(\k\m_j;\c^{j}_{\m}\bigr)
=(-1)^{\frac{\k\m_j}{2}+1-e(\c^{j}_{\m})}
\frac{\t(\c^{j}_{\m})}{2}\Bigl(\frac{2\pi}{N(\c^{j}_{\m})}\Bigr)^{\k\m_j}
\frac{B_{\k\m_j,\ol{\c^{j}_{\m}}}}{(\k\m_j)!}.
\end{equation}
 Substituting expression \eqref{for:L-middd} into equation \eqref{for:L},
 we obtain the desired formula.
\end{proof}

\subsection{Examples}

 For a partition $\m$, 
 we denote by
 $\prod^{\m,\odd}_{\m_j\ge a}$ (resp. $\prod^{\m,\even}_{\m_j\ge a}$)
 the product over all $1\le j\le \ell(\m)$
 such that $\m_j$ is odd (resp. even) and $\m_j\ge a$.
 In particular, we omit the condition``$\m_j\ge a$'' if $a=1$. 

\subsubsection{Principal Dirichlet characters}

 Let $\c^{(N)}_{1}$ be the principal Dirichlet character modulo $N$.
 Note that $\k=\k(k,\c^{(N)}_{1})=2k$.

\begin{cor}
\label{cor:principal}
 It holds that   
\begin{align}
\label{for:principal}
 L^{\om}_{d}\bigl(\{2k\}^d;\{\c^{(N)}_{1}\}^d\bigr)
=(-1)^{kd}2^{2kd}
\Biggl\{\sum_{\m\vdash d}
\frac{\e^{\om}_{\m}}{z_{\m}}\frac{(-1)^{\ell(\m)}}{2^{\ell(\m)}}
\Bigl(\prod_{p:\mathrm{prime}\atop p|N}
\prod^{\ell(\m)}_{j=1}\bigl(1-p^{-2k\m_j}\bigr)\Bigr)\wt{B}_{2k\m}
\Biggr\}\pi^{2kd}.
\end{align}
 In particular, we have 
 $L^{\om}_{d}(\{2k\}^d;\{\c^{(N)}_{1}\}^d)\in\bQ{\pi}^{2kd}$.
\end{cor}
\begin{proof}
 Since $(\c^{(N)}_1)^{j}_{\m}=\bsym{1}$ for all $1\le j\le \ell(\m)$,
 we have
 $N_{\m}(\c^{(N)}_1)=1$,
 $e_{\m}(\c^{(N)}_1)=0$,
 $\t_{\m}(\c^{(N)}_1)=1$ and
 $\wt{B}_{2k\m,\c^{(N)}_1}=\wt{B}_{2k\m}$.
 Hence the claim follows immediately 
 from expression \eqref{for:multiple Dirichlet}.    
\end{proof}
 
\begin{example}
[The case $\c=\c^{(1)}_1=\bsym{1}$]
 For $\k=2k$ with $k\in\bN$,
 we have 
\begin{equation}
\label{for:zeta}
 \z^{\om}_{d}\bigl(\{2k\}^d\bigr)
=L^{\om}_{d}\bigl(\{2k\}^d;\{\bsym{1}\}^d\bigr)
=(-1)^{kd}2^{2kd}
\Biggl\{
\sum_{\m\vdash d}
\frac{\e^{\om}_{\m}}{z_{\m}}\frac{(-1)^{\ell(\m)}}{2^{\ell(\m)}}\wt{B}_{2k\m}
\Biggr\}\pi^{2kd}.
\end{equation}
 The expression \eqref{for:zeta} for $\om=\star$
 gives the result obtained in \cite{Ohno}. 
\end{example}

\begin{example}
[The case $\c=\c^{(2)}_1$]
 For $\k=2k$ with $k\in\bN$,
 we have 
\begin{equation}
\label{for:c-2}
 L^{\om}_{d}\bigl(\{2k\}^d;\{\c^{(2)}_{1}\}^d\bigr)
=(-1)^{kd}\Biggl\{
\sum_{\m\vdash d}\frac{\e^{\om}_{\m}}{z_{\m}}
\frac{(-1)^{\ell(\m)}\prod^{\ell(\m)}_{j=1}\bigl(2^{2k\m_j}-1\bigr)}
{2^{\ell(\m)}}\wt{B}_{2k\m}
\Biggr\}\pi^{2kd}.
\end{equation}
\end{example}

\subsubsection{Primitive real Dirichlet characters}

 For a fundamental discriminant $D$,
 let $\c_{D}$ be the associated 
 primitive real Dirichlet character modulo $|D|$
 defined by the Kronecker symbol $(\frac{D}{\cdot})$
 (see \cite{Zagier1981}).
 Note that
 $e(\c_{D})=e(D):=\frac{1-\sgn{D}}{2}$
 where $\sgn{D}$ is the signature of $D$, whence 
 $\k=\k(k,\c_{D})=2k+e(D)$.

\begin{cor}
\label{cor:primitive-real}
 It holds that  
\begin{multline}
\label{for:primitive-real}
 L^{\om}_{d}\bigl(\{\k\}^d;\{\c_{D}\}^d\bigr)\\
=(-1)^{\frac{\k d}{2}}2^{\k d}\Biggl\{\sum_{\m\vdash
d}\frac{\e^{\om}_{\m}}{z_{\m}}\frac{(-1)^{\ell(\m)-\frac{1}{2}e(D)\ell(\m_{\mathrm{o}})}}{2^{\ell(\m)}}\frac{\prod_{p:\mathrm{prime}\atop
p||D|}\prod^{\m,\even}_{j}\bigl(1-p^{-\k\m_j}\bigr)}{|D|^{\k
 |\m_{\odd}|-\frac{1}{2}\ell(\m_{\odd})}}\wt{B}_{\k\m_{\odd},\c_D}\wt{B}_{\k\m_{\even}}\Biggr\}\pi^{\k d}.
\end{multline}
\end{cor}
\begin{proof}
 Since $\c_{D}$ is real, 
 we have $(\c_{D})^{j}_{\m}=\bsym{1}$ if $\m_j$ is even
 and $\c_{D}$ otherwise.
  One can therefore obtain the formula by using the well-known property
 $\t(\c_D)=i^{e(D)}|D|^{\frac{1}{2}}$. 
\end{proof}

\begin{example}
[The case $\c=\c_{-4}$]
 For $\k=2k+1$ with $k\in\bZ_{\ge 0}$,
 we have 
\begin{multline}
\label{for:c-4}
 L^{\om}_{d}\bigl(\{\k\}^d;\{\c_{-4}\}^d\bigr)
=(-1)^{\frac{\k d}{2}}
\Biggl\{\sum_{\m\vdash d}
\frac{\e^{\om}_{\m}}{z_{\m}}
\frac{(-1)^{\ell(\m)+\frac{1}{2}\ell(\m_{\odd})}
\prod^{\m,\even}_{j}\bigl(2^{\k\m_j}-1\bigr)}{2^{\ell(\m)+\k |\m_{\odd}|}}
\wt{E}_{\k\m_{\odd}-1}\wt{B}_{\k\m_{\even}}
\Biggr\}\pi^{\k d}.
\end{multline} 
\end{example}

\begin{remark}
 For a non-principal Dirichlet character $\c$,
 it is known that $L(1,\c)\ne 0$ (see, e.g., \cite{Zagier1981}).
 However, in general, it does not seem to be the case that
 $L^{\om}_{d}(\{1\}^d;\{\c\}^d)\ne 0$ for $d\ge 2$.
 In fact, when $\om=\bul$,
 we can obtain the following example of a pair of $d$ and $\c$
 such that $L^{\bul}_{d}(\{1\}^d;\{\c\}^d)=0$
 (although we were not able to derive a similar result for the  case $\om=\star$): 
\begin{align*}
 L^{\bul}_{2}\bigl(\{1\}^2;\{\c_{-8}\}^2\bigr)
&=\frac{1}{1\cdot 3}+\frac{-1}{1\cdot 5}+\frac{-1}{1\cdot
 7}+\frac{1}{1\cdot 9}+\frac{1}{1\cdot 11}+\frac{-1}{1\cdot
 13}+\frac{-1}{1\cdot 15}+\frac{1}{1\cdot 17}+\cdots\\
&\qquad\quad\ +\frac{-1}{3\cdot 5}+\frac{-1}{3\cdot
 7}+\frac{1}{3\cdot 9}+\frac{1}{3\cdot 11}+\frac{-1}{3\cdot
 13}+\frac{-1}{3\cdot 15}+\frac{1}{3\cdot 17}+\cdots\nonumber\\
&\qquad\qquad\qquad\ \,+\frac{1}{5\cdot 7}+\frac{-1}{5\cdot
 9}+\frac{-1}{5\cdot 11}+\frac{1}{5\cdot 13}+\frac{1}{5\cdot
 15}+\frac{-1}{5\cdot 17}+\cdots\nonumber\\[3pt]
&\qquad\qquad\qquad\qquad\quad\ \ +\cdots \nonumber\\
&=-\Bigl\{\frac{3}{4}\frac{B_{2}}{2!}-\frac{1}{16}\bigl(\frac{B_{1,\c_{-8}}}{1!}\bigr)^2\Bigr\}{\pi}^2=0.\nonumber
\end{align*}
 In the last equality,
 we have used the facts
 $B_2=\frac{1}{6}$ and $B_{1,\c_{-8}}=-1$.
\end{remark}

\subsection{$q$-Analogues of the multiple $L$-functions}

 Let $0<q<1$ and $[n]_q:=\frac{1-q^n}{1-q}$.
 Defined a $q$-analogue $L^{\om}_{q,d}$ of $L^{\om}_{d}$ by 
\[
 L^{\om}_{q,d}(s_1,\ldots,s_d;f_1,\ldots,f_d)
:=\sum_{(m_{1},\ldots,m_{d})\in I^{\om}_{d}}
\frac{f_{1}(m_{1})q^{m_1(s_1-1)}\cdots f_{d}(m_{d})q^{m_d(s_d-1)}}
{[m_{1}]_q^{s_{1}}\cdots [m_{d}]_q^{s_{d}}}.
\]
 Similarly to the case for $L^{\om}_{d}$ (or $q=1$),
 we write $L^{\om}_{q,d}$ as $\z^{\om}_{q,d}$ if $f_j=\bsym{1}$ for all $j$ 
 and $L_q(s;f)=L^{\om}_{q,1}(s;f)$.
 The function $L^{\om}_{q,d}$ is a natural extension
 of the $q$-analogue of the Riemann zeta function $\z_q(s):=L_q(s;\bsym{1})$
 studied in \cite{KanekoKurokawaWakayama2003}
 (see also \cite{KawagoeWakayamaYamasaki2007}).
 Many relations among the values $\z^{\om}_{q,d}$ at positive integers
 have been studied;
 see, e.g.,
 \cite{Bradley2005a,Bradley2005b,OkudaTakeyama,Zhao2007}
 for the case $\om=\bul$ and
 \cite{OhnoOkuda2007} for $\om=\star$.
 Specializing $x_n=\c(n)q^{n(s-1)}[n]^{-s}_q$ with $\Re(s)>1$
 in equation \eqref{for:ephp}
 and using  
\begin{align}
\label{for:Lq-symmetric}
 v^{\om}_{(d)}(\bsym{x})
=L^{\om}_{q,d}\bigl(\{s\}^d;\{\c\}^d\bigr)
\quad \textrm{and} \quad 
 p_{(r)}(\bsym{x})
=\sum^{d-1}_{l=0}\binom{d-1}{l}(1-q)^lL_q(ds-l;\c^d),
\end{align}
 we obtain the following expression
 which gives a $q$-analogue of expression \eqref{for:L}. 

\begin{prop}
 It holds that 
\begin{equation}
\label{for:qL}
 L^{\om}_{q,d}\bigl(\{s\}^d;\{\c\}^d\bigr)
=\sum_{\m\vdash
 d}\frac{\e^{\om}_{\m}}{z_{\m}}\prod^{\ell(\m)}_{j=1}\sum^{\m_j-1}_{l_j=0}\binom{\m_j-1}{l_j}(1-q)^{l_j}L_q(\m_js-l_j;\c^{\m_j}).
\end{equation}
\qed
\end{prop}

 It may be noted that
 a recursive expression of $\z^{\bul}_{q,d}(\{s\}^d)$
 was obtained in \cite[Theorem~$1$]{Bradley2005a}.

\begin{remark}
 Let
 $(P,Q)=(L^{\bul},L^{\star})$ or $(L^{\bul}_q,L^{\star}_q)$.
 Then, for $d\ge 1$, 
 it is straightforward to obtain from equation \eqref{for:gen} that
 $\sum_{a,b\ge 0\atop a+b=d}(-1)^aP_a(\{s\}^{a};\{\c\}^a)Q_b(\{s\}^{b};\{\c\}^b)=0$
 (here we understand that
 $(L^{\om})_c=L^{\om}_c$ and $(L^{\om}_q)_c=L^{\om}_{q,c}$).
 Further, one can obtain more explicit relation between
 $P_{a}(\{s\}^{a};\{\c\}^a)$ and $Q_{b}(\{s\}^{b};\{\c\}^b)$.
 In fact,
 for any
 $(P,Q)\in\{(L^{\bul},L^{\star}),(L^{\star},L^{\bul}),(L^{\bul}_{q},L^{\star}_{q}),(L^{\star}_{q},L^{\bul}_{q})\}$,
 it holds that 
\begin{equation}
\label{for:L-L^*}
 P_d\bigl(\{s\}^d;\{\c\}^d\bigr)
=\sum_{\m\vdash d}
\e_{\m}u_{\m}\prod^{\ell(\m)}_{j=1}
Q_{\m_j}\bigl(\{s\}^{\m_j};\{\c\}^{\m_j}\bigr).
\end{equation}
 This follows from the relations
 $e_{(d)}=\sum_{\m\vdash d}\e_{\m}u_{\m}h_{\m}$ and
 $h_{(d)}=\sum_{\m\vdash d}\e_{\m}u_{\m}e_{\m}$
 (see \cite[Example~$20$, p.\,$33$]{Macdonald1995})
together with relations \eqref{for:L-symmetric} and \eqref{for:Lq-symmetric}.
 Here, we put $u_{\m}:=\binom{\ell(\m)}{m_1(\m),\,m_2(\m),\ldots}$
 where $\binom{n}{a,\,b,\ldots}$ ($n=a+b+\cdots$) denotes
 the multinomial coefficient. 
\end{remark}

\section{Multiple $L$-values via generating functions}
\label{sec:infinite products}

\subsection{Product expressions of generating functions}

 The second method
 to evaluate the value $L^{\om}_{d}(\{\k\}^d;\{\c\}^d)$
 is well-known for the case $\c=\bsym{1}$.
 Namely, we start from the following identity
 which is obtained by specializing $x_n=\c(n)n^{-\k}$ and
 replacing $t$ by $t^{\k}$ in the generating functions \eqref{for:gen}:  
\begin{equation}
\label{for:L-generating}
 \sum^{\infty}_{d=0}
\e^{d}_{\om}L^{\om}_{d}\bigl(\{\k\}^d;\{\c\}^d\bigr)t^{\k d}
=\prod^{\infty}_{n=1}\Bigl(1-\frac{\c(n)t^{\k}}{n^{\k}}\Bigr)^{-\e_{\om}}.
\end{equation}
 The right-hand side of equation \eqref{for:L-generating}
 can then be written as follows. 

\begin{prop}
\label{prop:infinite product}
 It holds that
\begin{align}
\label{for:N=1}
 \prod^{\infty}_{n=1}\Bigl(1-\frac{\bsym{1}(n)t^{2k}}{n^{2k}}\Bigr)
&=\prod^{\infty}_{n=1}\Bigl(1-\frac{t^{2k}}{n^{2k}}\Bigr)
=\prod^{k}_{l=1}\frac{\sin{(\pi \z^{l}_{2k}t)}}{\pi \z^{l}_{2k}t},\\
\label{for:N=2}
 \prod^{\infty}_{n=1}\Bigl(1-\frac{\c^{(2)}_1(n)t^{2k}}{n^{2k}}\Bigr)
&=\prod^{\infty}_{n=1}\Bigl(1-\frac{t^{2k}}{(2n-1)^{2k}}\Bigr)
=\prod^{k}_{l=1}\cos\bigl(\frac{\pi \z^{l}_{2k}t}{2}\bigr)
\intertext{for $k\in\bN$ and, for a Dirichlet character $\c$ modulo $N\ge 3$,}
\label{for:N>=3}
 \prod^{\infty}_{n=1}\Bigl(1-\frac{\c(n)t^{\k}}{n^{\k}}\Bigr)
&=\bigl(N^{-1}2^{N-1}\bigr)^{\frac{1}{2}\k}\prod^{\bar{N}}_{j=1}\prod^{\k}_{l=1}\sin\frac{\pi\bigl(j-\c(j)^{\frac{1}{\k}}
 \z^{l}_{\k}t\bigr)}{N}.
\end{align}
 Here, $\z_{m}:=e^{\frac{2\pi i}{m}}$ for $m\in\bN$ and
 $\bar{N}:=\Gauss{\frac{N-1}{2}}$ with
 $\Gauss{x}$ being the largest integer not exceeding $x$.
 In equation \eqref{for:N>=3},
 we take the argument of $\c(j)$ satisfying $\c(j)\ne 0$
 to be $0\le \arg{\c(j)}<2\pi$. 
\end{prop}

 We need the following lemma
 (see \cite[Example~$11$, Chapter~$\mathrm{VI}$, p.\,$115$]{Bromwich1949}).  

\begin{lem}
\label{lem:infi-prod-formula}
 For $a_{i},b_{i}\in\bC$ $(1\le i\le l)$
 satisfying $\sum^{l}_{i=1}a_i=\sum^{l}_{i=1}b_i$,
 we have
\begin{equation}
\label{for:inf-prod}
 \prod^{\infty}_{n=1}\frac{(n+a_1)\cdots (n+a_l)}{(n+b_1)\cdots (n+b_l)}
=\frac{\G(1+b_1)\cdots \G(1+b_l)}{\G(1+a_1)\cdots \G(1+a_l)}.
\end{equation}
 \qed
\end{lem}

\begin{proof}
[Proof of Proposition~\ref{prop:infinite product}]
 Let $N\ge 3$. 
 By the periodicity of $\c$,
 the left-hand side of equation \eqref{for:N>=3} equals   
\begin{align}
\label{for:N>=3-mid}
 \prod^{\infty}_{n=1}
\prod^{N-1}_{j=1}\Bigl(1-\frac{\c(nN-(N-j))t^{\k}}{(nN-(N-j))^{\k}}\Bigr)
&=\prod^{\infty}_{n=1}
\prod^{N-1}_{j=1}\frac{\bigl(n-\frac{N-j}{N}\bigr)^{\k}-\c(j)\bigl(\frac{t}{N}\bigr)^{\k}}{\bigl(n-\frac{N-j}{N}\bigr)^{\k}}\\
&=\prod^{\infty}_{n=1}
\frac{\prod^{N-1}_{j=1}\prod^{\k}_{l=1}\bigl(n-\frac{N-j+\c(j)^{\frac{1}{\k}}\z^l_{\k}t}{N}\bigr)}
{\prod^{N-1}_{j=1}\bigl(n-\frac{N-j}{N}\bigr)^{\k}}\nonumber\\
&=\Bigl(\prod^{N-1}_{j=1}\G\bigl(\frac{j}{N}\bigr)\Bigr)^{\k}
\prod^{N-1}_{j=1}\prod^{\k}_{l=1}\G\Bigl(\frac{j-\c(j)^{\frac{1}{\k}}\z^l_{\k}t}{N}\Bigr)^{-1}.\nonumber
\end{align}
 In the last equality,
 we have used formula \eqref{for:inf-prod}
 (it is easy to check the condition in Lemma~\ref{lem:infi-prod-formula}
 from $\sum^{N-1}_{j=1}\c(j)=0$ when $\k=1$ and 
 $\sum^{\k}_{l=1}\z^{l}_{\k}=0$ otherwise).
 Note that  
 $\prod^{N-1}_{j=1}\G(\frac{j}{N})=(N^{-1}(2\pi)^{N-1})^{\frac{1}{2}}$
 by the Gauss-Legendre formula 
 $\G(Na)(2\pi)^{\frac{N-1}{2}}=N^{Na-\frac{1}{2}}\prod^{N-1}_{j=0}\G(a+\frac{j}{N})$
 with $a=1$.
 Further,
 the double product in the final right-hand side of \eqref{for:N>=3-mid}
 can be written as
\begin{align*}
 \prod^{\k}_{l=1}\G\Bigl(\frac{\frac{N}{2}-\c(\frac{N}{2})^{\frac{1}{\k}}\z^l_{\k}t}{N}\Bigr)^{-\d_{N}}
&\prod^{\bar{N}}_{j=1}\G\Bigl(\frac{j-\c(j)^{\frac{1}{\k}}\z^l_{\k}t}{N}\Bigr)^{-1}\G\Bigl(\frac{N-j-\c(N-j)^{\frac{1}{\k}}\z^{l}_{\k}t}{N}\Bigr)^{-1},
\end{align*}
 where $\d_{N}:=1$ if $N$ is even and $0$ otherwise.
 Notice also that
 $\c(\frac{N}{2})=0$ for even $N\ge 3$
 and  
 $\c(N-j)^{\frac{1}{\k}}\z^l_{\k}=-\c(j)^{\frac{1}{\k}}\z^{l-k}_{\k}$
 if $0\le \arg{\c(j)}<\pi$ and
 $-\c(j)^{\frac{1}{\k}}\z^{l+k}_{\k}$
 otherwise.
 Therefore,
 using the formula $\G(\frac{1}{2})=\sqrt{\pi}$ and
 replacing $l\pm k$ in $\z^{l\pm k}_{\k}$ by $l$,
 we see that the above expression is equal to
\[
 \pi^{-\frac{\d_N}{2}}
\prod^{\k}_{l=1}\prod^{\bar{N}}_{j=1}
\G\Bigl(\frac{j-\c(j)^{\frac{1}{\k}}\z^l_{\k}t}{N}\Bigr)^{-1}
\G\Bigl(1-\frac{j-\c(j)^{\frac{1}{\k}}\z^{l}_{\k}t}{N}\Bigr)^{-1}
=\pi^{-\frac{(N-1)}{2}\k}
\prod^{\k}_{l=1}\prod^{\bar{N}}_{j=1}
\sin\frac{\pi\bigl(j-\c(j)^{\frac{1}{\k}}
 \z^{l}_{\k}t\bigr)}{N}.
\]
 Here, we have used the reflection formula
 $\G(x)\G(1-x)=\frac{\pi}{\sin{(\pi x)}}$
 and the equality $-\frac{\d_N}{2}\k-\bar{N}\k=-\frac{(N-1)}{2}\k$.
 Hence,
 substituting the above expression into \eqref{for:N>=3-mid},
 we obtain the desired formula. 

 The equations \eqref{for:N=1} and \eqref{for:N=2}
 are easily obtained from the infinite product expressions  
 $\frac{\sin{(\pi t)}}{\pi t}=\prod^{\infty}_{n=1}(1-\frac{t^2}{n^2})$ and 
 $\cos{(\frac{\pi t}{2})}=\prod^{\infty}_{n=1}(1-\frac{t^2}{(2n-1)^2})$, 
 respectively
 (see \cite{Nakamura2005} for \eqref{for:N=1}).
 This completes the proof. 
\end{proof}

\subsection{Evaluation formula $\mathrm{II}$}

 Let $a\notin\pi\bZ$.
 For $\om\in\{\bul,\star\}$,
 we define the sequences $\{T^{\om}_n(a)\}_{n\ge 0}$
 by the expansions  
\begin{equation}
\label{for:series}
 \sin{(a+t)}
=\sum^{\infty}_{n=0}T^{\bul}_{n}(a)\frac{t^n}{n!}
\quad \textrm{and} \quad  
 \cosec{(a+t)}
=\sum^{\infty}_{n=0}T^{\star}_{n}(a)\frac{t^n}{n!}.
\end{equation}
 It is clear that 
 $T^{\bul}_{n}(a)=(-1)^{\frac{1}{2}n(n-1)}\tri_n{(a)}$
 where $\tri_n(a)=\sin{(a)}$ if $n$ is even and $\cos{(a)}$ otherwise.
 On the other hand,
 $T^{\star}_{n}(a)$ can be written as  
 $T^{\star}_{n}(a)=\cosec{(a)}\sum^{n}_{k=0}\binom{n}{k}i^{k}E_kh_{n-k}(\cot{(a)})$ 
 where $h_l(\a)$ is the polynomial in $\a$ of degree $l$ defined by
 $(1+\a\tan{(t)})^{-1}=\sum^{\infty}_{l=0}h_l(\a)\frac{t^l}{l!}$
 (note that $h_l$ is even if $l$ is and odd otherwise). 

 For a Dirichlet character $\c$ modulo $N$, define the sequence
 $\{A^{\om}_{n}(\k,\c)\}_{n\ge 0}$ as follows;

 $(\mathrm{i})$\ 
 The case $N=1$ (i.e., $\c=\bsym{1}$ and $\k=2k$)\,:\,
 For $n\in\bZ_{\ge 0}$,
 let $A^{\om}_{2n+1}(2k,\bsym{1})\equiv 0$ and    
\begin{align}
\label{for:A-1}
 A^{\om}_{2n}(2k,\bsym{1}):&=
\begin{cases}
 \DS{
\sum_{n_1,\ldots,n_{k}\ge 0\atop n_1+\cdots+n_{k}=n}
\binom{2n+k}{2n_1+1,\ldots,2n_{k}+1}\z^{\sum^{k}_{l=1}ln_l}_k}
 & \textrm{$(\om=\bul)$},\\
 \DS{
\sum_{n_1,\ldots,n_{k}\ge 0\atop n_1+\cdots+n_{k}=n}
\binom{2n}{2n_1,\ldots,2n_{k}}
\Bigl(\prod^{k}_{l=1}\bigl(2^{2n_l}-2\bigr)B_{2n_l}\Bigr) \z^{\sum^{k}_{l=1}ln_l}_k}
 & \textrm{$(\om=\star)$}.
\end{cases}
\end{align}

 $(\mathrm{ii})$\ 
 The case $N=2$ (i.e., $\c=\c^{(2)}_1$ and $\k=2k$)\,:\,
 For $n\in\bZ_{\ge 0}$,
 let $A^{\om}_{2n+1}(2k,\c^{(2)}_1)\equiv 0$ and 
\begin{align}
\label{for:A-2}
 A^{\om}_{2n}(2k,\c^{(2)}_{1}):&=
\begin{cases}
 \DS{
\sum_{n_1,\ldots,n_{k}\ge 0\atop n_1+\cdots+n_{k}=n}
\binom{2n}{2n_1,\ldots,2n_{k}}\z^{\sum^{k}_{l=1}ln_l}_k}
 & \textrm{$(\om=\bul)$},\\
 \DS{\sum_{n_1,\ldots,n_{k}\ge 0\atop n_1+\cdots+n_{k}=n}
\binom{2n}{2n_1,\ldots,2n_{k}}
\Bigl(\prod^{k}_{l=1}E_{2n_l}\Bigr)\z^{\sum^{k}_{l=1}ln_l}_k}
 & \textrm{$(\om=\star)$}.
\end{cases}
\end{align}

 $(\mathrm{iii})$\ 
 The case $N\ge 3$\,:\, Define
\begin{align}
\label{for:A-N-j}
 A^{\om}_{n}(j;\k,\c):
&=\sum_{n_1,\ldots,n_{\k}\ge 0\atop n_1+\cdots+n_{\k}=n}
\binom{n}{n_1,\ldots,n_{\k}}
\Bigl(\prod^{\k}_{l=1}T^{\om}_{n_l}\bigl(\frac{\pi j}{N}\bigr)\Bigr)
\z^{\sum^{\k}_{l=1}ln_l}_{\k}
 \quad \textrm{$(1\le j\le \bar{N})$}
\intertext{and let}
\label{for:A-N}
 A^{\om}_{n}(\k,\c):
&=\sum_{n_1,\ldots,n_{\bar{N}}\ge 0\atop n_1+\cdots+n_{\bar{N}}=n}
\binom{n}{n_1,\ldots,n_{\bar{N}}}
\prod^{\bar{N}}_{j=1}A^{\om}_{n_j}(j;\k,\c)\c(j)^{\frac{n_j}{\k}}.
\end{align}

\begin{example}
[The case $\c=\c_{-4}$]
 Let $E_n(x)$ be the Euler polynomial defined by
 $\frac{2e^{tx}}{e^t+1}=\sum^{\infty}_{n=0}E_n(x)\frac{t^n}{n!}$.
 Note that $E_n=2^nE_n(\frac{1}{2})$. 
 Then, since 
\[
 \cosec{\bigl(\frac{\pi}{4}+t\bigr)}
=\frac{\sqrt{2}}{\cos{(t)}+\sin{(t)}}
=\frac{\sqrt{2}(i+1)}{2}\frac{2e^{4it\cdot\frac{3}{4}}}{e^{4it}+1}
-\frac{\sqrt{2}(i-1)}{2}\frac{2e^{-4it\cdot\frac{3}{4}}}{e^{-4it}+1},
\]
 we see that 
 $T^{\star}_n(\frac{\pi}{4})=(-1)^{\frac{1}{2}n(n+1)}2^{\frac{4n+1}{2}}E_n(\frac{3}{4})$.
 While on the other hand, it is clear that 
 $T^{\bul}_n(\frac{\pi}{4})=2^{-\frac{1}{2}}(-1)^{\frac{1}{2}n(n-1)}$.
 Hence,
 from definitions \eqref{for:A-N-j} and \eqref{for:A-N},
 we have 
\begin{multline}
\label{for:A-4}
 A^{\om}_{n}(\k,\c_{-4})=A^{\om}_{n}(1;\k,\c_{-4})\\=
\begin{cases}
 \DS{
2^{-\frac{\k}{2}}\sum_{n_1,\ldots,n_{\k}\ge 0\atop n_1+\cdots+n_{\k}=n}
\binom{n}{n_1,\ldots,n_{\k}}(-1)^{\frac{1}{2}\sum^{\k}_{l=1}n_{l}(n_{l}-1)}
\z^{\sum^{\k}_{l=1}ln_{l}}_{\k}}
 & \textrm{$(\om=\bul)$},\\
 \DS{
2^{\frac{\k(2d+1)}{2}}\sum_{n_1,\ldots,n_{\k}\ge 0\atop n_1+\cdots+n_{\k}=n}
\binom{n}{n_1,\ldots,n_{\k}}\Bigl(\prod^{\k}_{l=1}E_{n_{l}}\bigl(\frac{3}{4}\bigr)\Bigr)
(-1)^{\frac{1}{2}\sum^{\k}_{l=1}n_{l}(n_{l}+1)}
 \z^{\sum^{\k}_{l=1}ln_{l}}_{\k}}
  & \textrm{$(\om=\star)$}.
\end{cases}
\end{multline}
\end{example}

 Put
 $\wt{\e}_{\om}:=\frac{\e_{\om}+1}{2}$,
 that is, $\wt{\e}_{\bul}=0$ and $\wt{\e}_{\star}=1$.
 We then obtain the following theorem 
 which gives extensions of the formulae for $\z^{\om}_d$ obtained in
 \cite{Nakamura2005,BorweinBradleyBroadhurst1997} for the case $\om=\bul$
 and \cite{Muneta} for $\om=\star$
 (see also \cite{Wenchang2000}). 

\begin{thm}
\label{thm:evaluation2}
 Let $\c$ be a Dirichlet character modulo $N$
 and $\k=2k+e(\c)\ge 1$. 
 Let
\[
 C^{\om}_{d}(\k,\c):=
\begin{cases}
 \DS{
\frac{\e^{d}_{\om}(-1)^{k(d-\wt{\e}_{\om})}}{((2d+1-\wt{\e}_{\om})k)!}}
 & \textrm{$(N=1,\ \k=2k)$},\\
 \DS{
\frac{\e^d_{\om}(-1)^{kd}}{2^{2kd}(2kd)!}}
 & \textrm{$(N=2,\ \k=2k)$},\\
 \DS{
\frac{\e^d_{\om}(-1)^{\k d}}{(N^{2d-\e_{\om}}2^{(N-1)\e_{\om}})^{\frac{\k}{2}}(\k d)!}}
 & \textrm{$(N\ge 3)$}.
\end{cases}
\]
 Then, we have $A^{\om}_{n}(\k,\c)=0$ if $\k\nmid n$ and 
\begin{equation}
\label{for:multiple Dirichlet 2}
 L^{\om}_d\bigl(\{\k\}^d;\{\c\}^d\bigr)
=C^{\om}_{d}(\k,\c)A^{\om}_{\k d}(\k,\c){\pi}^{\k d}.
\end{equation}
\end{thm}
\begin{proof}
 Let $N\ge 3$.
 It is straightforward
 from the definitions of $T^{\om}_n(a)$
 to see that
\begin{equation}
\label{for:sin-cosec}
 \prod^{\bar{N}}_{j=1}\prod^{\k}_{l=1}
\Bigl( \sin\frac{\pi\bigl(j-\c(j)^{\frac{1}{\k}}\z^{l}_{\k}t\bigr)}{N}
\Bigr)^{-\e_{\om}}
=\sum^{\infty}_{n=0}A^{\om}_{n}(\k,\c)\frac{(-1)^n}{n!}
\Bigl(\frac{\pi t}{N}\Bigr)^n.
\end{equation}
 In fact,
 when $\om=\bul$,
 the left-hand side of equation \eqref{for:sin-cosec} is given by  
\begin{align*}
 \prod^{\bar{N}}_{j=1}
\sum^{\infty}_{n_1,\ldots,n_{\k}=0}
\prod^{\k}_{l=1}T^{\bul}_{n_l}\bigl(\frac{\pi j}{N}\bigr)
\frac{1}{n_l!}\Bigl(-\frac{\pi
 \c(j)^{\frac{1}{\k}}\z^{l}_{\k}t}{N}\Bigr)^{n_l}
&=\prod^{\bar{N}}_{j=1}
\sum^{\infty}_{n=0}A^{\bul}_{n}(j;\k,\c)\frac{(-1)^n}{n!}
\Bigl(\frac{\pi \c(j)^{\frac{1}{\k}}t}{N}\bigr)^n\\
&=\sum_{n_1,\ldots,n_{\bar{N}}=0}
\prod^{\bar{N}}_{j=1}A^{\bul}_{n_j}(j;\k,\c)\frac{(-1)^n}{n!}
\Bigl(\frac{\pi \c(j)^{\frac{1}{\k}}t}{N}\bigr)^{n_j}\\
&=\sum^{\infty}_{n=0}A^{\bul}_{n}(\k,\c)\frac{(-1)^n}{n!}
\Bigl(\frac{\pi t}{N}\Bigr)^n.
\end{align*}
 The case when $\om=\star$ is similar.
 Therefore,
 from equation \eqref{for:N>=3}, we have 
\begin{equation}
\label{for:L-Taylor}
 \prod^{\infty}_{n=1}\Bigl(1-\frac{\c(n)t^{\k}}{n^{\k}}\Bigr)^{-\e_{\om}}
=\sum^{\infty}_{n=0}
\Biggl\{
(N^{-1}2^{N-1})^{-\frac{1}{2}\e_{\om}\k}
A^{\om}_{n}(\k,\c)\frac{(-1)^n}{n!}\Bigl(\frac{\pi }{N}\Bigr)^n
\Biggr\}t^n.
\end{equation}
 Hence,
 comparing the coefficients of $t^{\k d}$
 in expressions \eqref{for:L-generating} and \eqref{for:L-Taylor},
 we see that $A^{\om}_{n}(\k,\c)=0$ if $\k\nmid n$ and
 obtain formula \eqref{for:multiple Dirichlet 2} for $N\ge 3$.

 One can similarly obtain formulae \eqref{for:multiple Dirichlet 2}
 for the other cases, that is, $N=1$ or $N=2$, by using the expansions
\begin{align}
\label{for:many Taylor expansions}
\renewcommand{\arraycolsep}{1.5pt}
\begin{array}{rclcrcl}
 \sin{(t)}
&=&\DS{\sum^{\infty}_{n=0}\frac{(-1)^n}{(2n+1)!}t^{2n+1}},
 & \qquad &
 \cosec{(t)}
&=&\DS{\sum^{\infty}_{n=0}\frac{(-1)^{n-1}(2^{2n}-2)B_{2n}}{(2n)!}t^{2n-1}},\\
 \cos{(t)}
&=&\DS{\sum^{\infty}_{n=0}\frac{(-1)^n}{(2n)!}t^{2n}},
 & \qquad &
 \sec{(t)}
&=&\DS{\sum^{\infty}_{n=0}\frac{(-1)^{n}E_{2n}}{(2n)!}t^{2n}}.
\end{array}
\end{align}  
 Actually, for example let $N=1$.
 Then, from \eqref{for:N=1}, we have 
\begin{align}
\label{for:exN=1}
 \sum^{\infty}_{d=0}
\e^{d}_{\om}L^{\om}_{d}\bigl(\{2k\}^d;\{\bsym{1}\}^d\bigr)t^{2kd}
&=\prod^{\infty}_{n=1}\Bigl(1-\frac{\bsym{1}(n)t^{2k}}{n^{2k}}\Bigr)^{-\e_{\om}}
=\prod^{k}_{l=1}\Bigl(\frac{\sin{(\pi \z^{l}_{2k}t)}}{\pi \z^{l}_{2k}t}\Bigr)^{-\e_{\om}}\\
&=
\begin{cases}
  \DS{\sum^{\infty}_{n=0}A^{\bul}_{2n}(2k,\bsym{1})\frac{(-1)^n{\pi}^{2n}}{(2n+k)!}{t}^{2n}}
 & \textrm{$(\om=\bul)$},\\
 \DS{\sum^{\infty}_{n=0}A^{\star}_{2n}(2k,\bsym{1})\frac{(-1)^{n-k}{\pi}^{2n}}{(2n)!}{t}^{2n}} 
 & \textrm{$(\om=\star)$}.
\end{cases}
\nonumber
\end{align}
 In the last equality with $\om=\bul$ (resp. $\om=\star$),
 we have used the expansion of $\sin{(t)}$ (resp. $\cosec{(t)}$) in \eqref{for:many Taylor expansions}.
 Now comparing the coefficients of $t^{2kd}$ on both side of \eqref{for:exN=1},
 we have
\[
 L^{\om}_{d}\bigl(\{2k\}^d;\{\bsym{1}\}^d\bigr)\,\bigl(=\z^{\om}_d(\{2k\}^d)\bigr)
=
\begin{cases}
  \DS{\frac{\e^{d}_{\bul}(-1)^{kd}}{(2kd+k)!}A^{\bul}_{2kd}(2k,\bsym{1}){\pi}^{2kd}}
 & \textrm{$(\om=\bul)$},\\
 \DS{\frac{\e^{d}_{\star}(-1)^{kd-k}}{(2kd)!}A^{\star}_{2kd}(2k,\bsym{1}){\pi}^{2kd}} 
 & \textrm{$(\om=\star)$},
\end{cases}
\]
 whence the desired formula follows.
 This completes the proof of the theorem.
\end{proof}

\begin{example}
[The case $\c=\bsym{1}$]
 The following results are well known
 (see \cite{BorweinBradleyBroadhurst1997,Muneta}): 
\begin{center}
\begin{tabular}{c|c|c|c}
 $\k=2k$ & $2$ & $4$ & $6$ \\[1pt]
\hline
\hline
 $\z^{\bul}_d(\{\k\}^d)$
 & ${\frac{1}{(2d+1)!}\pi^{2d}}$
 & ${\frac{2^{2d+1}}{(4d+2)!}\pi^{4d}}$ 
 & ${\frac{3\cdot 2^{6d+1}}{(6d+3)!}\pi^{6d}}$
 \\[5pt]
\hline
 $\z^{\star}_d(\{\k\}^d)$  
 & ${\frac{(-1)^{d-1}(2^{2d}-2)B_{2d}}{(2d)!}\pi^{2d}}$ & $-$ & $-$ \\[5pt]
\end{tabular}
\end{center} 
 Notice that,
 from expression \eqref{for:multiple Dirichlet 2},
 we have
 $\z^{\star}_d(\{4\}^d)=((2^{4d}+4)S_{2d}(-1)-4S_{2d}(-4)){\pi}^{4d}/(4d)!$
 where $S_{k}(t):=\sum^{k}_{n=0}\binom{2k}{2n}t^{n}B_{2n}B_{2k-2n}$.
 However, we do not know whether
 the expression can be reduced to a ``simpler'' (or ``closed'') form,
 because it would appear to be difficult to simplify
 the sum $S_{k}(t)$ for a general $t\in\bC$
 (see \cite{Dilcher1996} for the comment on $S_{k}(-1)$).
 We remark that,  
 on the other hand,
 it is well-known that $S_{k}(1)=-(2k-1)B_{2k}$ for $k\ne 1$
 and $S_{1}(1)=\frac{1}{3}$.   
\end{example}

\begin{example}
[The case $\c=\c^{(2)}_{1}$]
 We obtain the following results:
\begin{center}
\begin{tabular}{c|c|c|c}
 $\k=2k$ & $2$ & $4$ & $6$ \\[1pt]
\hline
\hline
 $L^{\bul}_d(\{\k\}^d;\{\c^{(2)}_1\}^d)$
 & ${\frac{1}{2^{2d}(2d)!}\pi^{2d}}$
 & ${\frac{1}{2^{2d}(4d)!}\pi^{4d}}$ 
 & ${\frac{3}{4(6d)!}\pi^{6d}}$
 \\[5pt]
\hline
 $L^{\star}_d(\{\k\}^d;\{\c^{(2)}_1\}^d)$  
 & ${\frac{(-1)^{d}E_{2d}}{2^{2d}(2d)!}\pi^{2d}}$ & $-$ & $-$ \\[5pt]
\end{tabular}
\end{center}
 From expression \eqref{for:multiple Dirichlet 2},
 it can be expressed as     
 $L^{\star}_d(\{4\}^d;\{\c^{(2)}_1\}^d)=T_{2d}(-1)\pi^{4d}/(2^{4d}(4d)!)$
 where $T_{k}(t):=\sum^{k}_{n=0}\binom{2k}{2n}t^{n}E_{2n}E_{2k-2n}$. 
 We also note that 
 $T_{k}(1)=2^{2k+1}E_{2k+1}(1)$ for $k\ge 0$ (see \cite{Dilcher1996}).
\end{example}

\begin{example}
[The case $\c=\c_{-4}$]
 We obtain the following results:
\begin{center}
\begin{tabular}{c|c|c|c}
 $\k=2k+1$ & $1$ & $3$ & $5$ \\[1pt]
\hline
\hline
 $L^{\bul}_d(\{\k\}^d;\{\c_{-4}\}^d)$
 & ${\frac{(-1)^{\frac{1}{2}d(d-1)}}{2^{2d}d!}\pi^{d}}$
 & ${\frac{(-1)^{\frac{1}{2}d(d-1)}\cdot 3}{2^{3d+1}(3d)!}\pi^{3d}}$ 
 & ${\frac{(-1)^{\frac{1}{2}d(d-1)}\cdot 5(L_{5d}-1)}{2^{5d+2}(5d)!}\pi^{5d}}$
 \\[5pt]
\hline
 $L^{\star}_d(\{\k\}^d;\{\c_{-4}\}^d)$  
 & ${\frac{(-1)^{\frac{1}{2}d(d-1)}E_{d}(\frac{3}{4})}{d!}\pi^{d}}$ &
 $-$ & $-$ \\[5pt]
\end{tabular}
\end{center}
 Here,
 $L_{d}$ is the Lucas number defined by the recursion equation
 $L_{1}=1$, $L_{2}=3$ and $L_{d}=L_{d-1}+L_{d-2}$ for $d\ge 3$
 (it can be shown that $L_d=\a^d+\b^{d}$
 where $\a=\frac{1+\sqrt{5}}{2}$ and $\b=\frac{1-\sqrt{5}}{2}$).
 Note that it is obtained in \cite{BorweinBradleyBroadhurst1997} that 
 $\z^{\bul}_d(\{10\}^d)=2^{10d+1}\cdot 5(L_{10d+5}+1){\pi}^{10d}/(10d+5)!$, and 
 can be similarly shown that
 $L^{\bul}_d(\{10\}^d;\{\c^{(2)}_1\}^d)=5(L_{10d}+1){\pi}^{10d}/(2^4(10d)!)$. 
\end{example}

\section{Summation formulae for the Bernoulli and Euler numbers}
\label{sec:applications}

 From
 Theorem~\ref{thm:evaluation1} and Theorem~\ref{thm:evaluation2},
 we immediately obtain the following corollary:

\begin{cor}
 It holds that
\begin{equation}
\label{for:sum formula}
 \sum_{\m\vdash d}
\frac{\e^{\om}_{\m}}{z_{\m}}\frac{(-1)^{\ell(\m)-e_{\m}(\c)}}{2^{\ell(\m)}}
\frac{W_{\m}(\k;\c)\t_{\m}(\c)}{N_{\m}(\c)^{\k}}\wt{B}_{\k\m,\ol{\c}}
=(-1)^{\frac{\k d}{2}}2^{-\k d}C^{\om}_{d}(\k,\c)A^{\om}_{\k d}(\k,\c). 
\end{equation}
\qed
\end{cor} 
 
 The above equation expresses many (non-trivial) relations
 among the generalized Bernoulli numbers.
 We here give several formulae involving
 the Bernoulli and Euler numbers. 

\begin{example}
 Let $\c=\bsym{1}$ and $\c=\c^{(2)}_1$.
 Then, from expression \eqref{for:sum formula}, we have respectively 
\begin{align}
\label{for:sum-1}
 \sum_{\m\vdash d}\frac{\e^{\om}_{\m}}{z_{\m}}
\frac{(-1)^{\ell(\m)}}{2^{\ell(\m)}}\wt{B}_{2k\m}
&=\frac{\e^{d}_{\om}(-1)^{k\wt{\e}_{\om}}A^{\om}_{2kd}(2k,\bsym{1})}
{2^{2kd}((2d+1-\wt{\e}_{\om})k)!},\\
\label{for:sum-2}
 \sum_{\m\vdash d}\frac{\e^{\om}_{\m}}{z_{\m}}
\frac{(-1)^{\ell(\m)}\prod^{\ell(\m)}_{j=1}\bigl(2^{2k\m_j}-1\bigr)}{2^{\ell(\m)}}\wt{B}_{2k\m}
&=\frac{\e^d_{\om}A^{\om}_{2kd}(2k,\c^{(2)}_1)}
{2^{2kd}(2kd)!}.
\end{align}
 These expressions show that
 both finite sums
 $A^{\om}_{2kd}(2k,\bsym{1})$ and $A^{\om}_{2kd}(2k,\c^{(2)}_1)$
 can be written as sums of products of the Bernoulli numbers.
 On the other hand,
 putting $d=1$ and $\om=\bul$ in expressions \eqref{for:sum-1} and \eqref{for:sum-2},
 one can also see
 from definitions \eqref{for:A-1} and \eqref{for:A-2} that 
 $B_{2k}$ ($k\ge 1$) has the following expressions
 in terms of the multinomial coefficients and the $k$-th root of unity:  
\begin{align}
\label{for:sum-Bernoulli1}
 B_{2k}
&=\frac{(2k)!}{2^{2k-1}(3k)!}
\sum_{n_1,\ldots,n_{k}\ge 0\atop n_1+\cdots+n_{k}=k}
\binom{3k}{2n_1+1,\ldots,2n_{k}+1}\z^{\sum^{k}_{l=1}ln_l}_k,\\
\label{for:sum-Bernoulli2}
 B_{2k}
&=\frac{1}{2^{2k-1}(2^{2k}-1)}
\sum_{n_1,\ldots,n_{k}\ge 0\atop n_1+\cdots+n_{k}=k}
\binom{2k}{2n_1,\ldots,2n_{k}}\z^{\sum^{k}_{l=1}ln_l}_k.
\end{align} 
 Formula \eqref{for:sum-Bernoulli1}
 was essentially obtained by Nakamura \cite{Nakamura2005}
 (one can easily check that expression \eqref{for:sum-Bernoulli1}
 is equivalent to $(2.2)$ in \cite{Nakamura2005}).
 Furthermore, 
 putting $k=1$ and $\om=\star$
 in expression \eqref{for:sum-1} (resp. \eqref{for:sum-2})
 and noting
 $A^{\bul}_{2d}(2,\c^{(2)}_1)=(2^{2d}-2)B_{2d}$
 (resp. $A^{\star}_{2d}(2,\c^{(2)}_1)=E_{2d}$),
 we obtain an expression for $B_{2d}$ (resp. $E_{2d}$) ($d\ge 1$)
 in terms of a sum of products of the Bernoulli numbers:
\begin{align}
\label{for:sum-Bernoulli3}
 B_{2d}
&=-\frac{2^{2d}(2d)!}{(2^{2d}-2)}
\sum_{\m\vdash d}\frac{1}{z_{\m}}
\frac{(-1)^{\ell(\m)}}{2^{\ell(\m)}}
\prod^{\ell(\m)}_{j=1}\frac{B_{2\m_j}}{(2\m_j)!},\\
\label{for:sum-Euler1}
 E_{2d}
&=2^{2d}(2d)!
\sum_{\m\vdash d}\frac{1}{z_{\m}}
\frac{(-1)^{\ell(\m)}}{2^{\ell(\m)}}
\prod^{\ell(\m)}_{j=1}(2^{2\m_j}-1)
\prod^{\ell(\m)}_{j=1}\frac{B_{2\m_j}}{(2\m_j)!}.
\end{align}
 Formula \eqref{for:sum-Bernoulli3}
 was obtained by Ohno in \cite{Ohno2}.
\end{example}

%
%

\begin{example}
 Let $\c=\c_{-4}$.
 Then, formula \eqref{for:sum formula} yields 
\begin{equation}
\label{for:sum-4}
 \sum_{\m\vdash d}\frac{\e^{\om}_{\m}}{z_{\m}}
\frac{(-1)^{\ell(\m)+\frac{1}{2}\ell(\m_{\odd})}
\prod^{\m,\even}_{j}\bigl(2^{\k\m_j}-1\bigr)}{2^{\ell(\m)+\k |\m_{\odd}|}}
\wt{E}_{\k\m_{\odd}-1}\wt{B}_{\k\m_{\even}}
=\frac{\e_{\om}^{\k d}(-1)^{\frac{\k d}{2}}A^{\om}_{\k d}(\k,\c_{-4})}
{2^{2\k d(1-\wt{\e}_{\om})}(\k d)!}.
\end{equation}
 Similarly to expressions 
 \eqref{for:sum-Bernoulli1} and \eqref{for:sum-Bernoulli2},
 putting $d=1$ and $\om=\bul$
 and writing $\k=2k+1$, 
 we see from expression \eqref{for:A-4} that
 $E_{2k}$ ($k\ge 0$) can be written as 
\begin{equation}
\label{for:sum-Euler}
 E_{2k}
=\frac{(-1)^k}{(2k+1)2^{2k}}\sum_{n_1,\ldots,n_{2k+1}\ge 0\atop
 n_1+\cdots+n_{2k+1}=2k+1}\binom{2k+1}{n_1,\ldots,n_{2k+1}}(-1)^{\frac{1}{2}\sum^{2k+1}_{l=1}n_l(n_l-1)}\z^{\sum^{2k+1}_{l=1}ln_l}_{2k+1}.
\end{equation}
\end{example}

\section{Special values at non-positive integers}
\label{sec:non-positive}

 It has been shown that,
 in some special cases,
 $L^{\bul}_d(s_1,\ldots,s_d;\c_1,\ldots,\c_d)$
 can be continued meromorphically to the whole space $\bC^{d}$
 and has possible singularities
 on $s_{d}=1$ and
 $\sum^{j}_{k=1}s_{d-k+1}\in\bZ_{\le j}$ for $j=2,3,\ldots,d$
 (see \cite{AkiyamaEgamiTanigawa2001,Matsumoto2003,Zhao2000}
 when $\c_j=\bsym{1}$ for all $j$, and
 \cite{AkiyamaIshikawa2002}
 when $\c_1,\ldots,\c_d$ are of the same conductor.
 One can also obtain more precise information about the singularities
 from these references.).
 Hence, in such cases, it is easy to see that
 $L^{\star}_d(s_1,\ldots,s_d;\c_1,\ldots,\c_d)$
 also admits a meromorphic continuation to $\bC^{d}$
 with the same possible singularities
 since $L^{\star}_{d}$ can be expressed
 in terms of $L^{\bul}_{d'}$ for $1\le d'\le d$.   
 In this section,
 we study a special value,
 called a {\it central limit value},
 of $L^{\om}_d(s_1,\ldots,s_d;\{\c\}^d)$
 at $s_1=\cdots=s_d=-\k$
 where $\k=2k+e(\c)$ is a non-negative integer.
 Since it is seen that 
 such a point is a point of indeterminacy of the function,
 to define the value,
 we have to choose the limiting process
 of $(s_1,\ldots,s_d)\to (-\k,\ldots,-\k)$ 
 (for more precise details,
 see \cite{AkiyamaEgamiTanigawa2001,AkiyamaTanigawa2001}).
 Here,
 the central limit values
 $(L^{\om}_d)^{C}(s_1,\ldots,s_d;\c_1,\ldots,\c_d)$ are defined by 
\[
 (L^{\om}_d)^{C}\bigl(s_1,\ldots,s_d;\c_1,\ldots,\c_d\bigr)
:=\lim_{\d\to 0}L^{\om}_{d}\bigl(s_1+\d,\ldots,s_d+\d;\c_1,\ldots,\c_d\bigr).
\]
 We then obtain the following results,
 which give generalizations
 of the formulae for $\z^{\bul}_d$ 
 obtained in \cite[Corollary~$2$ \,(ii)]{Kamano2006}.

\begin{thm}
\label{thm:central values}
 $(\mathrm{i})$\ 
 For $N\ge 1$, it holds that
\begin{align}
\label{for:k=0}
 (L^{\om}_d)^{C}\bigl(\{0\}^d;\{\c^{(N)}_{1}\}^d\bigr)
&=\begin{cases}
 \DS{\e^d_{\om}(-1)^d\binom{\frac{\e_{\om}}{2}}{d}}
   & \textrm{$(N=1)$},\\
 0 & \textrm{$(N\ge 2)$},\\
\end{cases}\\ 
\label{for:-2k}
 (L^{\om}_d)^{C}\bigl(\{-2k\}^d;\{\c^{(N)}_{1}\}^d\bigr)
&=0 \quad (k\in\bN).
\end{align}
 $(\mathrm{ii})$\ 
 If $\c$ is non-principal,
 then,
 for $\k=2k+e(\c)\ge 0$,
 we have $(L^{\om}_d)^{C}(\{-\k\}^d;\{\c\}^d)=0$.  
\end{thm} 
\begin{proof}
 We start from the equation
$\sum^{\infty}_{d=0}\e^{d}_{\om}L^{\om}_d(\{s\}^d;\{\c\}^d)t^d=\exp(\e_{\om}\sum_{n\ge 1}\frac{1}{n}L(ns;\c^n)t^n)$,
 which is obtained 
 from expressions \eqref{for:gen} and \eqref{for:L-symmetric}.
 Note that this is valid for any $s\in\bC$
 unless $(s,\ldots,s)$ is a singularity of $L^{\om}_d$.
 Hence,
 using the formula $L(-\k,\c)=-\frac{B_{\k+1,\c}}{\k+1}$,
 we have
\begin{align}
\label{for:central}
 \sum^{\infty}_{d=0}
\e^{d}_{\om}(L^{\om}_d)^{C}\bigl(\{-\k\}^d;\{\c\}^d\bigr)t^d
=\exp\Bigl(-\e_{\om}
\sum^{\infty}_{n=1}\frac{B_{n\k+1,\c^n}}{n\bigl(n\k+1\bigr)}t^n\Bigr).
\end{align}

 We first assume that $\c$ is not principal.
 Then, noting $(n\k+1)+e(\c^n)=2n(k+e(\c))+1\equiv 1$ $\pmod{2}$
 (here we have used the identity $e(\c^n)\equiv n e(\c)$ $\pmod 2$),
 we have $B_{n\k+1,\c^n}=0$ for all $n\ge 1$
 because $B_{m,\c}=0$ if $m+e(\c)\equiv 1$ $\pmod 2$.
 This implies that 
 the right-hand side of equation \eqref{for:central}
 is identically equal to $e^{0}=1$,
 whence the claim $(\mathrm{ii})$ follows.

 We next let $\c=\c^{(N)}_1$. Note that  
\begin{equation}
\label{for:BerBer}
 B_{m,(\c^{(N)}_1)^n}=B_{m,\c^{(N)}_1}
=\prod_{p:\mathrm{prime}\atop p|N}(1-p^{m-1})\cdot B_{m}
\qquad (m\in\bZ_{\ge 0}).
\end{equation}
 Therefore,
 writing $\k=2k$ for $k\in\bN$ in equation \eqref{for:central} and
 using $B_{2nk+1}=0$,
 one similarly obtains equation \eqref{for:-2k}.
 To prove equation \eqref{for:k=0},
 we further set $\k=2k=0$ in equation \eqref{for:central}.
 When $N\ge 2$,
 since $B_{1,(\c^{(N)}_1)^n}=0$ by equation \eqref{for:BerBer},
 we obtain the claim.
 On the other hand when $N=1$,
 since $B_1=\frac{1}{2}$,
 the right-hand side of equation \eqref{for:central} can be written as 
\[
 \exp\Bigl(-\frac{\e_{\om}}{2}\sum_{n\ge 1}\frac{1}{n}t^n\Bigl)
=\exp\Bigl(\frac{\e_{\om}}{2}\log{(1-t)}\Bigl)
=(1-t)^{\frac{\e_{\om}}{2}}
=\sum^{\infty}_{d=0}(-1)^d\binom{\frac{e_{\om}}{2}}{d}t^d.
\]  
 This shows the claim.
\end{proof}

\begin{remark}
 By employing a ``renormalization procedures'' in quantum fields theory, 
 one can also define values of multiple zeta functions at non-positive integers. 
 For further details, see 
 \cite{{BerlineVergne2007},{GuoZhang2008a},{GuoZhang2008b},{ManchonPaycha}}
 and references therein.
\end{remark}

\section{Concluding remarks}
\label{sec:Concluding remarks}

\subsection{Multiple zeta functions attached to the Schur functions}

 It may be interesting to find
 a function $\z_{\l}(s_1,\ldots,s_d)$ for a partition $\l\vdash d$
 such that
 $\z_{(1^d)}=\z^{\bul}_d$ and $\z_{(d)}=\z^{\star}_d$.
 Recall that
 $\z^{\bul}_d(\{s\}^d)=e_d(\bsym{n}^{-s})$ and
 $\z^{\star}_d(\{s\}^d)=h_d(\bsym{n}^{-s})$
 where $\bsym{n}^{-s}=(1^{-s},2^{-s},\ldots)$.
 Therefore,
 since $s_{(1^d)}=e_d$ and $s_{(d)}=h_d$
 where $s_{\l}$ is the Schur function attached to $\l$
 (see \cite{Macdonald1995}),
 such a function $\z_{\l}$ can be regarded
 as a ``multiple zeta function attached to $s_{\l}$'':
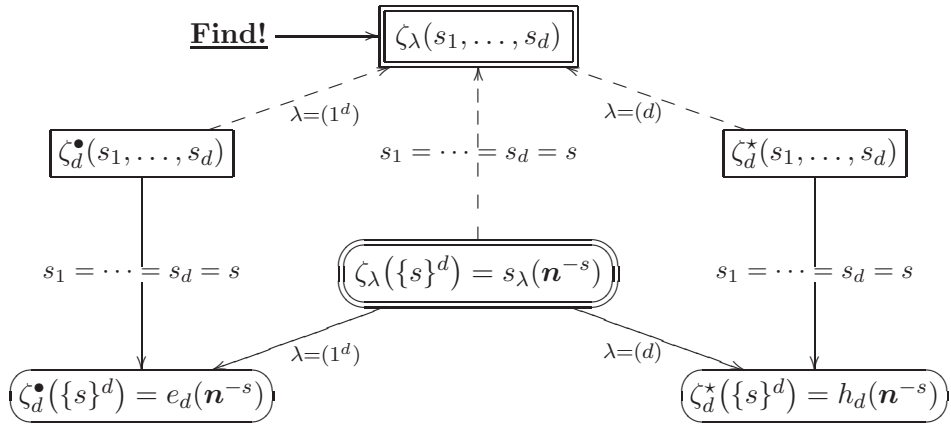
\begin{center}

\[
 \xymatrix{
 \qquad\qquad\qquad \ul{\txt{\bf Find!}} \ar[r] &
 *++[F=]{\z_{\l}(s_1,\ldots,s_d)} & \\
  *+[F-]{\z^{\bul}_{d}(s_1,\ldots,s_d)} \ar@{-->}[ru]_{\l=(1^d)}
 \ar[dd]|{\txt{{\small $s_1=\cdots=s_d=s$}}} & &
 *+[F-]{\z^{\star}_{d}(s_1,\ldots,s_d)} \ar@{-->}[lu]^{\l=(d)}
 \ar[dd]|{\txt{{\small $s_1=\cdots=s_d=s$}}} \\
 & *++[F=:<10pt>]{\z_{\l}\bigl(\{s\}^d\bigr)=s_{\l}(\bsym{n}^{-s})}
 \ar@{-->}[uu]|{\txt{{\small $s_1=\cdots=s_d=s$}}}
 \ar[ld]^{\l=(1^d)}
 \ar[rd]_{\l=(d)} & \\
 *+[F-:<10pt>]{\z^{\bul}_{d}\bigl(\{s\}^d\bigr)=e_d(\bsym{n}^{-s})} & &
 *+[F-:<10pt>]{\z^{\star}_{d}\bigl(\{s\}^d\bigr)=h_d(\bsym{n}^{-s})}
 }
\]
\end{center}

\begin{ques}
 For a fixed $\l$,
 are there any relations
 (such as a sum formula or a duality)
 among the values $\z_{\l}(k_1,\ldots,k_d)$,
 which interpolate the relations for $\z^{\om}_d(k_1,\ldots,k_d)$?    
\end{ques}

\subsection{Alternating multiple zeta functions}

 Let
\[
 \z^{\om}_d(\ol{s_1},\ldots,\ol{s_d})
:=\sum_{(m_1,\ldots,m_d)\in I^{\om}_d}
\frac{(-1)^{m_1}\cdots (-1)^{m_d}}{m_1^{s_1}\cdots m_d^{s_d}}
=L^{\om}_{d}\bigl(s_1,\ldots,s_d;\{\vp^{1}_{2}\}^d\bigr),  
\]
 where
 $\vp^{a}_{N}(m):=\z_N^{am}$ for $1\le a\le N$.
 This is an alternating analogue of $\z^{\om}_d(s_1,\ldots,s_d)$.
 When $\om=\bul$,
 the values at positive integers and their relations
 have been studied in \cite{BorweinBradleyBroadhurst1997,BowmanBradley2003}.
 Similarly to the proof of Theorem~\ref{thm:evaluation1},
 using expression \eqref{for:ephp}
 with the specialization $x_n=(-1)^nn^{-s}$,
 we get the following expressions of $\z^{\om}_d(\{\ol{s}\}^d)$
 in terms of a sum over partitions of $d$
 (note that this series also converges for $s=1$):
\begin{align*}
 \z^{\om}_d\bigl(\{\ol{1}\}^d\bigr)
&=\sum_{\m\vdash d}
\Biggl\{
\frac{\e^{\om}_{\m}}{z_{\m}}
\frac{(-1)^{\ell(\m)+\frac{|\m_{\even}|}{2}}\prod^{\m,\odd}_{\m_j\ge 3}
\bigl(2^{\m_j-1}-1\bigr)}{2^{\ell(\m_{\even})-\ell(\m_{\odd})-(|\m_{\even}|-|\m_{\odd}|)}}
\wt{B}_{\m_{\even}}\Biggr\}(\log{2})^{m_1(\m)}\cdot
\textstyle{\prod^{\m,\odd}_{\m_j\ge 3}}\z(\m_j)\cdot
 \pi^{|\m_{\even}|},\\
 \z^{\om}_d\bigl(\{\ol{2k}\}^d\bigr)
&=(-1)^{kd}2^{2kd}
\Biggl\{
\sum_{\m\vdash d}
\frac{\e^{\om}_{\m}}{z_{\m}}
\frac{(-1)^{\ell(\m_{\even})}\prod^{\m,\odd}_{j}
\bigl(2^{2k\m_j-1}-1\bigr)}{2^{\ell(\m_{\even})+2k|\m_{\odd}|}}
\wt{B}_{2k\m}\Biggr\}\pi^{2kd}\in\bQ{\pi}^{2kd}
\quad \textrm{$(k\in\bN)$}.
\end{align*}
 These follow from the identities $L(1;\vp^{1}_2)=-\log{2}$ and 
\begin{align*}
 L\bigl(2k\m_j;(\vp^{1}_2)^{\m_j}\bigr)=
\begin{cases}
 \DS{-\frac{2^{2k\m_j-1}-1}{2^{2k\m_j-1}}\z(2k\m_j)}
 & \textrm{if $\m_j$ is odd},\\
 \z(2k\m_j)
 & \textrm{otherwise}.
\end{cases}
\end{align*}

 Furthermore, since
 $\prod^{\infty}_{n=1}(1-\frac{\vp^{1}_{2}(n)t^{2k}}{n^{2k}})^{-\e_{\om}}
=\prod^{\infty}_{n=1}(1-\frac{(\z_{4k}t)^{2k}}{(2n-1)^{2k}})^{-\e_{\om}}
\prod^{\infty}_{n=1}(1-\frac{1}{n^{2k}}(\frac{t}{2})^{2k})^{-\e_{\om}}$,
 then, upon using formulae
 \eqref{for:N=1}, \eqref{for:N=2} and \eqref{for:many Taylor expansions},
 we have
\begin{equation}
\label{for:alternative}
 \z^{\om}_d\bigl(\{\ol{2k}\}^d\bigr)
=\frac{\e^d_{\om}(-1)^{k(d-\wt{\e}_{\om})}
A^{\om}_{2kd}(2k,\vp^{1}_{2})}{2^{2kd}((2d+1-\wt{\e}_{\om})k)!}{\pi}^{2kd},
\end{equation} 
 where $\{A^{\om}_{n}(2k,\vp^{1}_{2})\}_{n\ge 0}$ is defined by
 $A^{\om}_{2n+1}(2k,\vp^{1}_{2})\equiv 0$ and  
\[
 A^{\om}_{2n}(2k,\vp^{1}_{2}):
=\sum_{p,q\ge 0\atop p+q=n}\binom{2n+(1-\wt{\e}_{\om})k}{2p,2q+(1-\wt{\e}_{\om})k}\z^{p}_{2k}A^{\om}_{2p}(2k,\c^{(2)}_{1})A^{\om}_{2q}(2k,\bsym{1}).
\]  
  Formula \eqref{for:alternative} can be straightforwardly obtained from $(35)$
 in \cite{BorweinBradleyBroadhurst1997}.

\subsection{A higher-rank generalization}

 Let $s_{ij}\in\bC$ and $f_{ij}$ be arithmetic functions
 for $1\le i\le r$ and $1\le j\le d$.
 Write $S=(s_{ij})$ and $F=(f_{ij})$ ($r\times d$ matrices).
 We introduce the following ``higher-rank'' multiple $L$-function: 
\begin{align*}
 {}_rL^{\om}_d\bigl(S;F\bigr)
:=\sum_{(m^{1},\ldots,m^{d})\in{}_rI^{\om}_{d}}
\prod^{d}_{j=1}\frac{f_{1j}(m_{1j})\cdots f_{rj}(m_{rj})}
{m_{1j}^{s_{1j}}\cdots m_{rj}^{s_{rj}}},
\end{align*}
 where
 $m^{j}=(m_{1j},\ldots,m_{rj})\in\bN^{r}$ for $1\le j\le d$ and
 ${}_rI^{\om}_{d}\subset (\bN^{r})^d$ is defined by 
\[
 {}_rI^{\om}_d
:=
\begin{cases}
 \bigl\{(m^1,\ldots,m^d)\in(\bN^{r})^{d}\,|\,m^1<\cdots<m^d\bigr\} &
 \ \ \textrm{$(\om=\bul)$},\\
 \bigl\{(m^1,\ldots,m^d)\in(\bN^{r})^{d}\,|\,m^1\le\cdots\le m^d\bigr\} &
 \ \ \textrm{$(\om=\star)$}.
\end{cases}
\]
 Here, we are employing the lexicographic order on $\bN^r$.
 Namely,
 $(m_1,\ldots,m_r)>(n_1,\ldots,n_r)$
 if $m_1>n_1$, or
 $m_1=n_1$ and $m_2>n_2$, or
 $m_1=n_1$, $m_2=n_2$ and $m_3>n_3$, and so on.
 Clearly,
 ${}_1I^{\om}_d=I^{\om}_d$ and ${}_1L^{\om}_d=L^{\om}_d$.
 As for the case $r=1$,
 we write ${}_r\z^{\om}_d(S)$ as ${}_rL^{\om}_d(S;(\bsym{1}))$. 
 

 By induction on $d$, 
 it is shown that the function ${}_rL^{\om}_d(S;F)$
 can be expressed as a polynomial
 in $L^{\om}_{d'}(t_1,\ldots,t_{d'};g_1,\ldots,g_{d'})$ for $1\le d'\le d$
 where
 $t_l$ is a sum of $s_{ij}$ and
 $g_l$ is a product of $f_{ij}$.
 For example, we see that
\begin{align*}
 {}_rL^{\om}_1\bigl((s_{i1});(f_{i1})\bigr)
&=\prod^{r}_{i=1}L(s_{i1};{f_{i1}}),\\
{}_rL^{\om}_2\bigl((s_{ij});(f_{ij})\bigr)
&=\sideset{}{'}\sum^{r}_{p=1}\prod^{p-1}_{i=1}
L(s_{i1}+s_{i2};f_{i1}f_{i2})
L^{\bul}_2(s_{p1},s_{p2};f_{p1},f_{p2})
\prod^{r}_{i=p+1}L(s_{i1};f_{i1})L(s_{i2};f_{i2}).
\end{align*}
 Here,
 the prime $'$ means that we replace $\bul$ by $\om$
 in the summand for $p=r$.
 Hence,
 if $f_{ij}$ are bounded for all $i,j$,
 then we see
 from the expression of ${}_rL^{\om}_d$ in terms of $L^{\om}_{d'}$
 that 
 the series ${}_rL^{\om}_d(S;F)$ converges absolutely
 for $\Re(s_{1j})\ge 1$ for $1\le j\le d-1$ and $\Re(s_{ij})>1$ otherwise.

\begin{ques}
 Are there any relations
 among the values ${}_r\z^{\om}_{d}((k_{ij}))$,
 which are generalizations of the relations for $r=1$?    
\end{ques}

 Let us denote by $\{x_1,\ldots,x_r\}^d$
 the $r\times d$ matrix $(x_{ij})$ such that $x_{i1}=\cdots=x_{id}=x_i$
 for $1\le i\le r$.
 Let $\c_{i}$ be a Dirichlet character modulo $N_i$ and
 $\k_i=\k_{i}(k_i,\c)=2k_i+e(\c_i)$ for $1\le i\le r$.
 Though the function ${}_{r}L^{\om}_{d}$ is not new in the above sense,
 one can calculate the special values
 ${}_{r}L^{\om}_{d}(\{\k_1,\ldots,\k_r\}^d;\{\c_1,\ldots,\c_r\}^d)$
 by our first method.
 In fact, we finally can obtain the following proposition:

%

\begin{prop}
 Let $|\k|:=\k_1+\cdots+\k_r$.
 Then, we have 
\begin{multline}
\label{for:MLV-higher}
\ \ {}_{r}L^{\om}_{d}
\bigl(\{\k_1,\ldots,\k_r\}^d;\{\c_1,\ldots,\c_r\}^d\bigr)\\
=(-1)^{\frac{|\k| d}{2}}2^{|\k| d}
\Biggl\{
\sum_{\m\vdash d}\frac{\e^{\om}_{\m}}{z_{\m}}
\frac{(-1)^{r\ell(\m)-\sum^{r}_{i=1}e_{\m}(\c_i)}}{2^{r\ell(\m)}}
\prod^{r}_{i=1}\frac{W_{\m}(\k_i;\c_i)\t_{\m}(\c_i)}{N_{\m}(\c_i)^{\k_i}}
\wt{B}_{\k_i\m,\ol{\c_i}}\Biggr\}\pi^{|\k| d}. \ \ 
\end{multline}
 In particular,
 for $k_1,\ldots,k_r\in\bN$ with $|k|:=k_1+\cdots+k_r$,
 it holds that 
\[
 {}_{r}\z^{\om}_{d}
\bigl(\{2k_1,\ldots,2k_r\}^d\bigr)
=(-1)^{|k|d}2^{2|k|d}
\Biggl\{
\sum_{\m\vdash d}\frac{\e^{\om}_{\m}}{z_{\m}}
\frac{(-1)^{r\ell(\m)}}{2^{r\ell(\m)}}\prod^{r}_{i=1}
\wt{B}_{2k_i\m}\Biggr\}\pi^{2|k|d}
\in\bQ{\pi}^{2|k|d}. 
\]
\end{prop}
\begin{proof}
 Replacing the variable
 $x_n$ by $x^{1}_{n_1}\cdots x^{r}_{n_r}$ in \eqref{for:gen},
 we have 
\begin{equation}
\label{for:higher-gen}
 \sum^{\infty}_{d=0}\e_{\om}^{d}v^{\om}_{(d)}(\bsym{y})t^d
=\prod^{\infty}_{n_1,\ldots,n_r=1}\bigl(1-x^{1}_{n_1}\cdots x^{r}_{n_r}t\bigr)^{-\e_{\om}}
=\sum^{\infty}_{d=0}\e_{\om}^d\Bigl\{\sum_{\m\vdash
 d}\frac{\e^{\om}_{\m}}{z_{\m}}\prod^{r}_{i=1}p_{\m}(\bsym{x}^i)\Bigr\}t^d, 
\end{equation}
 where
 $\bsym{y}=(x^{1}_{n_1}\cdots x^{r}_{n_r})_{n_1,\ldots,n_r\ge 1}$
 (we also arrange the variables in $\bsym{y}$
 in the lexicographic order with respect to $(n_1,\ldots,n_r)$)
 and
 $\bsym{x}^i=(x^{i}_{n_i})_{n_i\ge 1}$.
 Then,
 specializing $x^{i}_{n_i}=\c_i(n_i)n_i^{-\k_i}$ in \eqref{for:higher-gen},
 we see that the left-hand side of equation \eqref{for:MLV-higher}
 is equal to  
 $\sum_{\m\vdash d}\frac{\e^{\om}_{\m}}{z_{\m}}\prod^{r}_{i=1}\prod^{\ell(\m)}_{j=1}L(\k_i\m_j;\c_i^{\m_j})$.
 Therefore,
 one can obtain equation \eqref{for:MLV-higher}
 by following the same approach as
 in the proof of Theorem~\ref{thm:evaluation1}. 
\end{proof}

\begin{remark}
 From expression \eqref{for:MLV-higher},
 we see that  
\[
 {}_{r}L^{\om}_{d}
\bigl(\{\k_{\s(1)},\ldots,\k_{\s(r)}\}^d;\{\c_{\s(1)},\ldots,\c_{\s(r)}\}^d\bigr)
={}_{r}L^{\om}_{d}
\bigl(\{\k_{1},\ldots,\k_{r}\}^d;\{\c_{1},\ldots,\c_{r}\}^d\bigr)
\]
 for any $\s\in\fS_r$
 where $\fS_r$ is the symmetric group of degree $r$. 
\end{remark}

\begin{remark}
 It appears to be difficult to evaluate the values
 ${}_{r}L^{\om}_{d}(\{\k_1,\ldots,\k_r\}^d;\{\c_1,\ldots,\c_r\}^d)$ 
 for $r\ge 2$
 by using the second method of the present study
 (that is, by using the expansions \eqref{for:many Taylor expansions}
 of trigonometric functions). This is because the
  generating function obtained from \eqref{for:higher-gen}
for these values is a ``multiple'' infinite product.

\end{remark}

\begin{Acknowledgement}
 The author would like to thank Professors Sho Matsumoto and Kazufumi Kimoto for valuable comments.
 He is also thankful to the referee for his/her careful reading and useful comments. 
\end{Acknowledgement}


\bigskip

\noindent
\textsc{Yoshinori YAMASAKI}\\
 Faculty of Mathematics, Kyushu University.\\
 Hakozaki, Fukuoka, 812-8581 JAPAN.\\
 \texttt{yamasaki@math.kyushu-u.ac.jp}\\

\noindent
{\it Current address:}\\
 Graduate School of Science and Engineering, Ehime University.\\
 Bunkyo-cho, Matsuyama, 790-8577 JAPAN.\\
 \texttt{yamasaki@math.sci.ehime-u.ac.jp}

\end{document}